\documentclass[12pt, twoside, leqno]{article}

\usepackage{amsmath,amsthm}
\usepackage{amssymb}
\usepackage{amsrefs}
\usepackage{enumerate}

\usepackage{graphicx}

\newtheorem{theorem}{Theorem} 
\newtheorem{corollary}{Corollary}

\newtheorem{exam}{Example}

\newtheorem*{rem}{Remark}

\newcommand\commentout[1]{}

\newcommand{\length}{\operatorname{length}}

\begin{document}
\title{On cubic Thue equations and the common index divisors of cyclic cubic fields}

\author{Mohammed Seddik\\
Universit\'e d'\'Evry Val d'Essonne, Universit\'e Paris-Saclay\\
Laboratoire de MathŽmatiques et Mod\'elisation d'\'Evry (UMR 8071)\\
I.B.G.B.I., 23 Bd. de France, 91037 \'Evry Cedex, France\\
 \hskip2cm mohammed.seddik@univ-evry.fr}

\date{}

\maketitle


\renewcommand{\thefootnote}{}



\renewcommand{\thefootnote}{\arabic{footnote}}
\setcounter{footnote}{0}


\begin{abstract}
In this paper, we investigate the common index divisors of cyclic cubic fields.   Let $a,b,c,d$ and $k$ are  integers,   we then solve the following  Thue cubic  equations:
\[ax^3+bx^2y+cxy^2+dy^3=  k\  \]
when $a,bc+d$ are odd and $3$ doesn't divide $v_2(k)$. 
\end{abstract}

\section{\textbf {Introduction and statement of main result}}
In 1909, Thue \cite{Thue} proved the following important result : If $f(x,y)\in \mathbb{Z}[x,y]$ is an irreducible binary form of degree at least 3, and $k$ non-zero integer, then the Diophantine equation $f(x,y)=k$
has only a finite number of integer solutions $(x,y)$. However, the method used does not help in finding the values of $x$ and $y$ in question. It is only possible to obtain an upper bound for the number of solutions. Baker \cite{Baker} has shown that bounds can also be given for the magnitude of the solutions. The bounds for $\arrowvert x\arrowvert$ and $\arrowvert y\arrowvert$ are extremely big, even for small values of degree of $f(x,y)$.\\

If $f(x,y)$ is an irreducible binary cubic form with negative discriminant, Delauney \cite{delone} and Nagell \cite{nagell} showed that the equation $f(x,y)=1$ has at most five integer solutions $(x,y)$.
Now if its discriminant is positive, then Evertse \cite{evertse}  showed that the equation $f(x,y)=1$ has at most twelve integer solutions $(x,y)$. Recently, Bennett \cite{Benett} refined  Delauney-Nagell-Evertse result as follows: if $f(x,1)$ has at least two distinct complex roots, then the equation $f(x,y)=1$ possesses at most $10$ solutions in integers $x$ and $y$.\\

Again, let $f(x,y)$ be an irreducible binary cubic form, the general equation $f(x,y)=k$ with $\vert D(f)\vert >\gamma k^{33}$,  $k$ non-negative integer, and $\gamma$ certain positive constant, Siegel proved there is at most $18$ solutions if $D(f)$ is positive (resp. Fjellstedt  proved there is at most $14$ solutions  when $D(f)$ is negative). See \cite[p.208]{nagell}.\\ 

In 1984, Ayad \cite{Ayad} proved that if  $f(x,y)$ is a binary form of degree $3$ with coefficients in $\mathbb{Z}$, $ Aut(f)$ its automorphisms group  and $H(f)$ its Hessian, then $Aut(f)$ is trivial except when $H(f)=\lambda g(x,y)$, $\lambda\in \mathbb{Z}^*$ and $g(x,y)$ is equivalent to $x^2-xy-y^2$. In this last case, $Aut(f)$ is cyclic of order $3$ and $f$ is equivalent to one binary form of type :
$$f_{m,n}(x,y)=mx^3-nx^2y-(n+3m)xy^2-my^3,\quad m,n\in\mathbb{Z},$$
so, the number of representations of integer $k$ by $f(x,y)$ is divisible by $3$. Note that the case $m=k=1$  is proved by Avanesov \cite{Avanesov}.\\

Recently, in \cite{Wakabayashi} Wakabayashi  proved that for $k=1$, the Thue equation $f_{m,n}(x,y)=1$ 
has at most three integer solutions except for a few known cases. Using the Pad\'e approximation method he obtained a lower bound for the size of solutions.\\

However, only few results are known when the group of automorphisms is trivial.
The family of cubic Thue equations $f_{1,n}(x,y)=\pm 1$, with $n\geq -1$, was studied by Thomas in \cite{Thomas}. Using Baker's method, Thomas proved that it has only trivial solutions except for a finite number of values of the parameter $n$, explicitly for $n<10^8$. Later, Mignotte solved the remaing cases for  this family of equations in \cite{Mignotte}.\\

Let $n$ be a rational integer and $\mathbb{K}_n=\mathbb{Q}(\theta)$ be a cyclic cubic number field generated by a root $\theta$ of $f_{1,n}(X,1)=X^3-nX^2-(n+3)X-1$ and let $\mathbb{O}_{\mathbb{K}_{n}}$ be its ring of integers. The polynomial $f_{1,n}(X,1)$ has discriminant $(n^2+3n+9)^2$. If $n^2+3n+9$ is square-free, then we have the discriminant of $\mathbb{K}_n$, $D(\mathbb{K}_n)=(n^2+3n+9)^2$ and $\mathbb{O}_{\mathbb{K}_n}=\mathbb{Z}[\theta]$ (there exist infinitely many such $n$, cf. Cusick \cite[pp. 63-73]{Cusick}).\\

  For $\theta$ a root of $P_n$, P\"etho and Lemmermeyer in \cite{Petho2} proved that for all $\alpha\in\mathbb{Z}[\theta]$ either $|N(\alpha)|\geq 2n+3$, or $\alpha$ is associated to an integer. Moreover, if $|N(\alpha)|=2n+3$, then $\alpha$ is associated to one of the conjugates of $\alpha-1$.  In \cite{Petho}, Mignotte, P\"etho and Lemmermeyer, by Baker's method and used the results of \cite{Petho2} and \cite[pp.63-73]{Cusick} solved completely the case $m=1$, $n\geq -1$ and $1\leq k\leq 2n+3$. Lettl, P\"etho and Voutier \cite{Lettl}, following an idea of Chudnovsky \cite{Chudnovsky}, improved the usual estimate of the Pad\'e approximation and gave a good upper bound for the size of solutions for $f_{1,n}(x,y)=k$, $n\geq 30$ and $k$ is arbitrary.  Wakabayashi \cite{Wakabayashi1} studied Thue inequality $|f_{m,n}(x,y)|\leq k$ with two parameters $m,n$.\\
  
In 2011, A. Hoshi \cite{Hoshi} studied the case when $k$ is a positive divisor of $n^2 + 3n + 9$,  and gave  a correspondence between integer solutions to the parametric family of cubic Thue equations
\[x^3-nx^2y-(n+3)xy^2-y^3=k\]
and isomorphism classes of the simplest cubic fields.  For more details on the study of simplest cubic fields see \cite{Shanks}.\\

A. Togb\'e \cite{Togbe}, using Baker's method and the results obtained by L. C. Washington \cite{Washington} and O. Lecacheux \cite{Lecacheux},  solved the family of parametrized Thue equations 
\[
x^3-(n^3-2n^2+3n-3)x^2y-n^2xy^2-y^3=\pm 1,\quad\text{when}\;n\geq 1.
\]
A. Togb\'e \cite{TogbŽ} using Baker's method and the results obtained by Y. Kishi \cite{Kishi}, solved the family of parametrized Thue equations 
\[ x^3-n(n^2+n+3)(n^2+2)x^2y-(n^3+2n^2+3n+3)xy^2-y^3=\pm 1, \quad
\text{when}\;n\geq 0.\]

Recently, Balady \cite{Balady}  showed that there are many families of cubics polynomials irreducibles with square discriminants.  We consider the family of polynomials \[X^3 +(n^7 +2n^6 +3n^5-n^4 -3n^3 -3n^2+ 3n +3)X^2+(-n^4 +3n)X-1\]
used in \cite{Balady} which satisfy the hypothesis of our main result.\\

Our method of proof is different to the previous methods, we use the common index divisors of cyclic cubic fields. Now we state our main result. 

\begin{theorem}\label{thm1}
 Let $a,b,c,d$ and $k$ be integers where $3\nmid v_2(k)$ and $a,bc+d$ are odd. Suppose that the polynomial $aX^3+bX^2+cX+d$ is irreducible with square discriminant. Then the cubic Thue equations
\[ ax^3+bx^2y+cxy^2+dy^3=  k\]
has no integer solution $(x,y)$.

\end{theorem}
Below we mention some immediate consequences of our result on the resolution of a parametric families of cubic equations.

\begin{corollary}\label{coro1}
 Let $k\in\mathbb{Z}$ where $3\nmid v_2(k)$. Then the parametric families of cubic Thue equations
\begin{enumerate}
\item $mx^3-nx^2y-(n+3m)xy^2-my^3=k, \;m,n\in\mathbb{Z},\;m\;\text{odd},\;\gcd(m,n)=1,$
\item $x^3-(n^3-2n^2+3n-3)x^2y-n^2xy^2-y^3=k,\;\;n\in\mathbb{Z},\;   n\not\equiv1\;mod\;4,$
\item $x^3-n(n^2+n+3)(n^2+2)x^2y-(n^3+2n^2+3n+3)xy^2-y^3=k,\;n\in\mathbb{Z},$
\item $x^3+(n^7 +2n^6 +3n^5-n^4 -3n^3 -3n^2+ 3n +3)x^2y+(-n^4 +3n)xy^2-y^3=k,\;n\in\mathbb{Z},$
\end{enumerate}
have no integer solution $(x,y)$.

\end{corollary}

\begin{corollary}\label{coro2}
 Let $a,b,c,d$ and $k$ be integers where $3\nmid v_2(k)$ and $a,bc+d$ are odd. Suppose that the polynomial $aX^3+bX^2+cX+d$ is irreducible with square discriminant. Then the homogeneous form
$$ ax^3+bx^2y+cxy^2+dy^3=kz^3$$
has only integer solution $(x,y,z)=(0,0,0)$.
\end{corollary} 

\begin{corollary}\label{coro3}
Let $a,b,c,d$ and $t$ be integers where $v_2(t)=1$ and $a,bc+d$ are odd. Suppose that the polynomial $aX^3+bX^2+cX+d+t$ is irreducible with square discriminant. Then the equation
\begin{eqnarray}\label{cubic}
ax^3+bx^2+cx+d= ty^2
\end{eqnarray}
has no integer solution $(x,y)$.
\end{corollary}
The  equation \eqref{cubic} was studied by Mordell \cite[pp.255-261]{Mordell}. He proved the following important result: if $a,b,c,d, t$ are arbitrary integers and the polynomial $ax^3+bx^2+cx+d$ has no squared linear  factor in $x$,  then the equation \eqref{cubic} has only a {\it finite number of integer solutions.} Here in under the by hypothesis of Corollary \ref{coro3}, the equation has no integer solution.

\end{document}